\documentclass[11pt]{article}
\usepackage[utf8]{inputenc}
\usepackage[T1]{fontenc}
\usepackage{amsmath}
\usepackage{amsfonts}
\usepackage{amssymb}
\usepackage{graphicx}

\usepackage{authblk}

\begin{document}

\title{\bf Nontrivial single axiom schemata and their quasi-nontriviality of Le\'{s}niewski-Ishimoto's propositional ontology $\bf L_1$}

\author[1]{\Large Takao Inou\'{e} \thanks{Corresponding Author: takaoapple@gmail.com}}
\author[2]{\Large Tadayoshi Miwa \thanks{miwa.tadayoshi@mail.u-tokyo.ac.jp}}
\affil[1]{\large Faculty of Informatics, Yamato University, Osaka, Japan\footnote{inoue.takao@yamato-u.ac.jp; (Personal) takaoapple@gmail.com [I prefer my personal mail.]}}
\affil[2]{\large Library, The University of Tokyo, Tokyo, Japan}

\date{February 2, 2025 (The 6th version)}


\maketitle

\abstract{On March 8, 1995, was found  the following \it nontrivial \rm single axiom-schema characteristic of Le\'{s}niewski-Ishimoto's propositional ontology $\bf L_1$ (Inou\'{e}, 1995b \cite{inoue16}).
$$(\mathrm{A_{M8})} \enspace \epsilon ab \wedge \epsilon cd . \supset . \epsilon aa \wedge \epsilon cc \wedge (\epsilon bc \supset . \epsilon ad \wedge \epsilon ba).$$

\medskip

In this paper, we shall present the progress about the above axiom-schema from 1995. Here we shall give two criteria \it nontiriviality \rm and \it quasi-nontriviality \rm in order to distinguish two axiom schemata. As main results, among others,  in \S 6 - \S 8, we shall give the simplified axiom schemata  ($\rm A_{S1}$), ($\rm A_{S2}$), ($\rm A_{S3N}$) and ($\rm A_{S3Nd}$) based on ($\mathrm{A_{M8}}$), their nontriviality and quasi-nontriviality. In \S 9 - \S 11, we shall give a lot of conjectures for nontrivial single axiom schemata for $\bf L_1$. We shall conclude this paper with summary and some remarks in \S 12.}

\medskip

\noindent SMC2020: 03B20, 03B60, 03F03, 03A05, 03A99
	
\tableofcontents

 \normalsize
	
\newtheorem{theorem}{Theorem}[section]
\newtheorem{proposition}{Proposition}[section]
\newtheorem{definition}{Definition}[section]
\newtheorem{corollary}{Corollary}[section]
\newtheorem{lemma}{Lemma}[section]
\newtheorem{convention}{convention}[section]
\newtheorem{remark}{Remark}[section]


\section{Introduction} 

On March 8, 1995, was found the following \it nontrivial \rm single axiom-schema characteristic of Le\'{s}niewski-Ishimoto's propositional ontology $\bf L_1$ (Inou\'{e}, 1995b \cite{inoue16}).

\bigskip

($\mathrm{A_{M8}}$) $\enspace$ $\epsilon ab \wedge \epsilon cd . \supset . \epsilon aa \wedge \epsilon cc \wedge (\epsilon bc \supset . \epsilon ad \wedge \epsilon ba)$.

\bigskip

In this paper, we shall present the progress about the above axiom-schema from 1995.

We shall propose natural criterions in order to distinguish one single axiom-schema with another and we shall see natural relationships among the proposed axiom schemata by means of the criterion. That is, here we shall give two criteria \it nontiriviality \rm and \it quasi-nontriviality \rm in order to distinguish two axiom schemata.

In the rest of this introduction, we shall give some preliminaries for $\mathbf{L_1}$. In \S 2, we shall provide the definition of nontrivial axiom schemata. In \S 3, the criterion of quasi-nontriviality is given in order to distinguish two nontrivial axiom schemata for $\mathbf{L_1}$. In \S 4, we shall give certain relations and properties about such criteria, that is, their conceptual structures. In \S 5, we shall present a $\mathbf{L_1}$-provable formula to deduce transitivity (Ax2 below) and exchangeability (Ax3 below) of axiom schemata for $\mathbf{L_1}$. In \S 6 and \S 7, we shall give three simplified axiom schemata ($\rm A_{S1}$), ($\rm A_{S2}$), ($\rm A_{S3N}$) and ($\rm A_{S3Nd}$) based on ($\mathrm{A_{M8}}$) and their nontriviality. In \S 8, as the main results of this paper, we shall show their quasi-nontriviality among them. In \S 9 - \S 11, we shall give a lot of conjectures for nontrivial single axiom schemata for $\bf L_1$. We shall conclude this paper with summary and  some remarks in \S 12.

Let us recall a formulation of $\mathbf{L_1}$, which was introduced in Ishimoto 1977 \cite{ishimoto1}. $\mathbf{L_1}$ is a propositional subsystem of  Le\'{s}niewski's ontology $\mathbf{L}$ (for $\mathbf{L}$, see Indrzejczak 2022 \cite{Indrzejczak2022}, Iwanu\'{s} 1972 \cite{iwa},  Lejewski 1958 \cite{lejewski}, S\l upecki 1954 \cite{slu}, Smirnov 1986 \cite{smirnov1}, Srzednicki \& Rickey 1984 \cite{sr}, Stachniak 1981 \cite{stach}, Surma \& Srzednicki 1992 \cite{surma}, Urbaniak 2014 \cite{urbaniak-book} and so on). 
The language of $\bf L_1$ consists of an infinite list of name variables $a, b, c, \dots$, $x_1, x_2, \dots$, $y_1, y_2, \dots$, Le\'{s}niewski's epsilon $\epsilon$, primitive logical symbols $\vee$ (disjunction) and $\neg$ (negation) and a set of auxiliary symbols \{ $)$, $($ \}. We shall first define formulas of $\bf L_1$.
We shall often use  the same infinite list of name variables as meta-name variables ranging the name variables of the list.
\medskip

\begin{definition} \rm The set of formulas of $\bf L_1$ is the smallest set $X$ which satisfies the following properties:
	
	$(1)$ For every pair of name variables $a$ and $b$, $\epsilon ab \in X$.
	
	$(2)$ If $A \in X$, then $(\neg A) \in X$.
	
	$(3)$ If $A, B \in X$, then $(A \vee B) \in X$.
	
\end{definition}
\medskip

Usually parentheses are omitted unless ambiguity arises. Other logical symbols are defined in terms of $\vee$ and $\neg$ as usual. So every atomic formula of $\bf L_1$ is of the form $\epsilon ab$ for some name variables $a$ and $b$, where $\epsilon$ is Le\'{s}niewski's epsilon. A very informal interpretation of $\epsilon ab$ in English may be the following:
$$\mbox{The $a$ is $b$}$$
(for this, see e.g. Prior 1965 \cite{prior65}).\footnote{Although this interpretation will not be used in this paper, it formed our intuitive ground to understand $\epsilon$. }

A Hilbert-style system for $\bf L_1$ is defined as the smallest set of formulas generated by 

\bigskip

\noindent $\mathsf{Axiom-schemata:}$

\medskip

(Ax1) \enspace $\epsilon ab$ $\supset$ $\epsilon aa$

\medskip

(Ax2) \enspace $\epsilon ab$ $\wedge$ $\epsilon bc$. $\supset$ $\epsilon ac$

\medskip

(Ax3) \enspace $\epsilon ab$ $\wedge$ $\epsilon bc$. $\supset$ $\epsilon ba$

\bigskip 

\noindent $\mathsf{Rule}$ :

\medskip

(R1) \enspace $\vdash_{\bf L_{1}}$ $A$, $\vdash_{\bf L_{1}}$ $A \supset B$ $\Rightarrow$ $\vdash_{\bf L_{1}}$ $B$,

\bigskip

\noindent where $a$, $b$, $c$ are meta-variables ranging name variables, and all the instances of classical tautology being closed under modus ponens. The axiom-schema (Ax3) can be replaced by the following simplified one due to Kanai 1989 \cite{kanai}:

\medskip 

(Ax3s) \enspace $\epsilon ab \wedge \epsilon bb. \supset \epsilon ba$,

\medskip

\noindent which is in some case advantageous.
\medskip

\noindent \bf Convention\rm . We shall use the following notation for a uniform subsitution for a formula (schemata) $A$,
$$\sigma = \begin{pmatrix}
u_1 & u_2 & \cdots & u_n \\
v_1 & v_2 & \cdots & v_n \\
\end{pmatrix}, \quad \sigma(A)$$
if $\sigma$ is a uniform substitution for variables (formulas) $u_1, u_2, \dots, u_n$ such that $u_i$ in $A$ is uniformly and simultaneously replaced by $v_i$ for each integer $1 \leq i \leq n$, and the result of the substitution is $\sigma(A)$.
\medskip

\begin{proposition}\label{Proposition 1}
For any uniform substitution $\sigma$ for meta-name variables, we have
$$\vdash_{\bf L_{1}} A \Rightarrow \enspace \vdash_{\bf L_{1}} \sigma(A).$$
\end{proposition}
Proof. Induction on derivation. $\Box$
\medskip

\begin{proposition}
For any multiple uniform substitution $\sigma$ for meta-name variables, we have
$$\vdash_{\bf L_{1}} A \Rightarrow \enspace \vdash_{\bf L_{1}} \sigma(A).$$
\end{proposition}
Proof. From Proposition \ref{Proposition 1}. $\Box$
\medskip

For the other topics for $\bf L_1$, e.g.  the interpretation of Le\'{s}niewski's $\epsilon$, axiomatic rejection, model theory, tableau systems and the modal interpretation of $\bf L_1$, some applications to linguistics and related studies, the reader will be recommended to refer to Blass 1994 \cite{blass1}, Inou\'{e} 1995a, 1995c, 1995d, 2021a, 2021b \cite{inoue15, inoue17, inoue1, inoue2021, inoue-Blass}, Inou\'{e} et al. 2021 \cite{Inoue-Ishi-Koba}, Ishimoto 1986, 1997 \cite{ishimoto4, Ishimoto1997}, Kobayashi \& Ishimoto 1982 \cite{koishi}, Ozawa \& Waragai 1985 \cite{OzawaWragai1985}, Pietruszczak1991  \cite{Piet1991}, Smirnov 1986 \cite{smirnov1}, Stachniak 1981 \cite{stach}, Takano 1985 \cite{takano0} and so on.


\section{The definition of nontrivial axiom schemata}

We shall, in this section, propose a natural criterion to distinguish a single axiom-schema with the set of  original axiom schemata, that is \it nontriviality. \rm 

\medskip 

\begin{definition} \rm Let $A$ be a single axiom-schema of $\mathbf{L_1}$. By $nv(A)$, 
we denote the ordered tuple of meta-variable of name variables occurring in $A$, 
where the order is based on the first occurence of the variable from the left.
\end{definition}
\medskip

For example, $nv(\epsilon ab \wedge \epsilon bc. \supset \epsilon ac) = (a, b, c).$
\medskip

\begin{definition} \rm Let $A$ be a single axiom-schema of $\mathbf{L_1}$. By $\# nv(A)$, 
we denote the length of the ordered tuple $nv(A)$.
\end{definition}
\medskip

For example, $\# nv(\epsilon ab \wedge \epsilon bc. \supset \epsilon ac) = 3.$
\medskip

Let us see the following single axiom-schema characteristic of $\bf L_1$.
\bigskip  

($\rm A_t$) $\enspace$ $\epsilon ab \supset . \epsilon aa \wedge (\epsilon bc \supset . \epsilon ac \wedge \epsilon ba)$.

\bigskip 

\begin{proposition} \rm 
$(A_t) \equiv . (Ax1) \wedge (Ax2) \wedge (Ax3)$ \it is an instance of tautology.
\end{proposition}
\medskip
Proof.  We easily see 
$$(P \supset . Q \wedge (R \supset . S \wedge T)) \equiv (P \supset Q . \wedge (P \wedge R . \supset S) \wedge (P \wedge R . \supset T))$$
is an instance of tautology, where $P, Q, R, S, T$ are distinct propositional variables. 
Take the following uniform substitution to the above formula:
$$\sigma_t = \begin{pmatrix}
P & Q & R & S & T \\
\epsilon ab & \epsilon aa & \epsilon bc & \epsilon ac & \epsilon ba \\
\end{pmatrix} . \enspace $$
\noindent $\Box$

We need some definitions for our further understanding of single axiom schemata characteristic of $\bf L_1$. We shall first in particular repeat to note our convention as follows: 

$$\sigma = \begin{pmatrix}
x_1 & x_2 & \cdots & x_n \\
y_1 & y_2 & \cdots & y_n \\
\end{pmatrix}$$
if $\sigma$ is a uniform substitution for (meta) name variables $x_1, x_2, \dots, x_n$ such that $x_i$ is replaced by $y_i$ for each integer $1 \leq i \leq n$
\medskip.

\begin{definition}\label{Definition 4} \rm Let $A$ be a single axiom-schema characteristic of $\mathbf{L_1}$ such that  $\# nv(A) \geq 3$. Say 
$$nv(A) = (x_1, x_2, x_3, x_4, \dots, x_n),$$
$$nv(\mathrm{(A_t)}) = nv(\epsilon ab \supset . \epsilon aa \wedge (\epsilon bc \supset . \epsilon ac \wedge \epsilon ba)) = (a, b, c),$$ and 
$$\# nv(\mathrm{(A_t)}) = 3,$$
where $n \geq 3$ and $x_1, x_2, x_3, x_4, \dots, x_n$ are mutually different. $A$ is \it trivial \rm (\it with respect to \rm $\mathrm{(A_t)}$\rm ) if 
there is a uniform substitution $\sigma$ to $A$ such that for some permutation $\rho$ of $(1, 2, 3, 4, \dots, n)$ and some set of 
mutually different meta-name variables $\{y_1, y_2, \dots, y_m \}$  ($m \geq 0$) with 
$$ \{x_1, x_2, x_3, x_4, \dots, x_n \} \cap \{y_1, y_2, \dots, y_m \} = \emptyset,$$
$$ \{a, b, c \} \cap \{y_1, y_2, \dots, y_m \} = \emptyset,$$
$$\sigma = \begin{pmatrix}
x_{\rho(1)} & x_{\rho(2)} & x_{\rho(3)} & x_{\rho(4)} & \dots & x_{\rho(n)} \\
a & b & c & y_1 & \dots & y_m \\
\end{pmatrix}, $$
$$\sigma(A) \equiv \mathrm{(A_t)}$$
is an instance of tautology of classical propositional logic, where $\equiv$ is the logical symbol for equivalence. If $A$ is not trivial (with respect to $\mathrm{(A_t)}$), it is said to be \it nontrivial \rm (\it with respect to $\mathrm{(A_t)}$\rm). (The end of \bf Definition \ref{Definition 4}\rm )
\end{definition}
\medskip

Informally, a trivial single axiom-schema characteristic of $\mathbf{L_1}$ is obtained by propositionally equivalent transformations from 
$$\mathrm{(Ax1)} \wedge \mathrm{(Ax2)} \wedge \mathrm{(Ax3)}$$ after some suitable uniform substitution for meta-name variables. So we are not interested in trivial single axiom schemata characteristic of $\mathbf{L_1}$. Nontrivial single axiom schemata are different from a set of original ones and are meaningful.

On March 8, 1995, was found the following \it nontrivial \rm single axiom-schema characteristic of $\bf L_1$ (Inou\'{e} 1995b \cite{inoue16}).

\bigskip

($\mathrm{A_{M8}}$) $\enspace$ $\epsilon ab \wedge \epsilon cd . \supset . \epsilon aa \wedge \epsilon cc \wedge (\epsilon bc \supset . \epsilon ad \wedge \epsilon ba)$.

\bigskip

\bf Remark 2.1\rm . Why does this Definition \ref{Definition 4} need such procedures? Probably, the reader will think of it so. Please think that for example, $\mathrm{A_{M8}}$ is an axiom-schema with respect to name variables occured in the schema. So for instance, the following schema $\sigma$($\mathrm{A_{M8}}$),
that is,
$$\epsilon ac \wedge \epsilon bd . \supset . \epsilon aa \wedge \epsilon bb \wedge (\epsilon cb \supset . \epsilon ad \wedge \epsilon ca).$$
\noindent is also ($\mathrm{A_{M8}}$) with
$$\sigma = \begin{pmatrix}
a & b & c & d \\
a & c & b & d \\
\end{pmatrix}.$$

\noindent (The end of \bf Remark 2.1\rm )
\medskip

\begin{proposition}\label{Proposition 4} 
$\mathrm{(A_{M8})}$  is nontrivial \rm (\it with respect to $\mathrm{(A_t)}$\rm).
\end{proposition}
\medskip

Proof.\footnote{There are many alternative proofs. This is a possible one.}
We shall show that is nontrivial. We see 
$$nv(\mathrm{A_{M8}}) = (a, b, c, d), \enspace \# nv(\mathrm{A_{M8}})= 4.$$ 

In order to prove the nontriviality of ($\rm A_{M8}$), we must check 24 ($= 4!$) cases, that is, 
the number of permutations of $(1, 2, 3, 4)$. However, it is quite easy if we observe the form of ($\rm A_{M8}$) and ($\rm A_{t}$).  We use $t$ (= true) and $f$ (= false) as usual.
The idea is to give a sentential valuation $v$ such that, noticing the form as $\epsilon \alpha \alpha$, the premise of $\sigma(\rm A_{M8})$ has the value $t$ and the conclusion of $\sigma(\rm A_{M8})$ has the value $f$ with respect to 
each permutation $\sigma$ of meta-name variables, whereas $v(\mathrm{(A_t)}) = t$ holds under the $v$.

Take a meta-name variable $y$ such that $\{a, b, c, d\} \cap \{y\} = \emptyset$.

\medskip

(Case 1) Let
$$\sigma = \begin{pmatrix}
u & v & c & w \\
a & b & c & y \\
\end{pmatrix}.$$

Then take a sentential valuation $v$ such that $v(\epsilon cc) = f$. For the rest $t$ is assigned.

Then we easily see $v((\rm A_t)) = \it t$ and $v(\sigma(\rm A_{M8})) = \it f$ under the assignment $v$. Thus $\sigma(\rm A_{M8}) \equiv (\rm A_t)$ is not an instance of tautology. 
\medskip


(Case 2) Let $x \neq c$.\footnote{In this case, we may choose a simpler presentation. However we shall take a present proof for the reader's intelligibility.}

\noindent (Subcase 2.1) Let 
$$\sigma = \begin{pmatrix}
c & v & x & w \\
a & b & c & y \\
\end{pmatrix}.$$

Then take a sentential valuation $v$ such that $v(\epsilon bb) = v(\epsilon cc) = v(\epsilon yy) = f$. For the rest $t$ is assigned.

Then we easily see $v((\rm A_t)) = \it t$ and $v(\sigma(\rm A_{M8})) = \it f$ under the assignment $v$. So $\sigma(\rm A_{M8}) \equiv (\rm A_t)$ is not an instance of tautology. 



\noindent (Subcase 2.2) Let 
$$\sigma = \begin{pmatrix}
u & c & x & w \\
a & b & c & y \\
\end{pmatrix}.$$

Then take a sentential valuation $v$ such that $v(\epsilon bb) = f$. For the rest $t$ is assigned.

Then we easily see $v((\rm A_t)) = \it t$ and $v(\sigma(\rm A_{M8})) = \it f$ under the assignment $v$. Hence $\sigma(\rm A_{M8}) \equiv (\rm A_t)$ is not an instance of tautology. 



\noindent (Subcase 2.3) Let 
$$\sigma = \begin{pmatrix}
u & v & x & c \\
a & b & c & y \\
\end{pmatrix}.$$

Then take a sentential valuation $v$ such that $v(\epsilon yy) = f$. For the rest is assigned to $t$.

Then we easily see $v((\rm A_t)) = \it t$ and $v(\sigma(\rm A_{M8})) = \it f$ under the assignment $v$, Thus $\sigma(\rm A_{M8}) \equiv (\rm A_t)$ is not an instance of tautology.  $\Box$
\medskip

For the reader who does not immediately see that ($\rm A_{M8}$) is a single axiom-schema characteristic of $\bf L_1$, we shall verify it below. 
Take the following uniform substitutions:

$$\sigma_1 = \begin{pmatrix}
a & b & c & d \\
a & b & a & b \\
\end{pmatrix},$$

$$\sigma_2 = \begin{pmatrix}
a & b & c & d \\
a & b & b & c \\
\end{pmatrix}.$$

By propositional logic, we immediately derive (Ax1) from $\sigma_1(\rm A_{M8})$. By the obtained (Ax1) and propositional logic, we immediately derive (Ax2) and (Ax3) from $\sigma_2(\rm A_{M8})$. On the other hand, $(\rm A_{M8})$ is a theorem of $\bf L_1$. Indeed, it is proved by the tableau method (for the method, see Kobayashi \& Ishimoto 1982 \cite{koishi} or Inou\'{e} 1995a \cite{inoue15} or Inou\'{e} et al. 2021 \cite{Inoue-Ishi-Koba}), or we can directly derive ($\rm A_{M8}$) from (Ax1)--(Ax3). Let us carry out the latter now. By (Ax1), we have 

\medskip

(C1a) $\epsilon ab \supset \epsilon aa$.

\medskip

(C1c) $\epsilon cd \supset \epsilon cc$.

\medskip

\noindent By (Ax2), (Ax3) and propositional logic, we easily get

\medskip

(C2) $\epsilon ab \wedge \epsilon bc \wedge \epsilon cd. \supset \epsilon ad$, 

\medskip

(C3) $\epsilon ab \wedge \epsilon bc \wedge \epsilon cd. \supset \epsilon ba$.

\medskip

\noindent Then we can derive $ (\rm A_{M8})$ from (C1a), (C1c), (C2) and (C3). 

\section{The criterion, the quasi-nontiriviality}


We shall give a criterion, that is, quasi-nontriviality (more general nontriviality than Definition \ref{Definition 4}) 
in order to distinguish one nontrivial single axiom-schema with another. The idea of the definition of quasi-nontriviality is 
that two single axiom schemata are compared on the ground of the same meta-variables.
\medskip

\begin{definition}\label{Definition 5} \rm Let $A$ and $B$ be single axiom schemata characteristic of $\mathbf{L_1}$ which is nontrivial (with respect to $\mathrm{(A_t)}$). Say 
$$nv(A) = (x_1, x_2, x_3, x_4, \dots, x_n),$$
$$nv(B) = (y_1, y_2, y_3, y_4, \dots, y_m),$$ where
$n \geq 3$, $m \geq 3$. $A$ is \it quasi-trivial \rm (\it with respect to $B$\rm ) if 

(Case 1) When $n \leq m$, 
there is a uniform subsititution $\sigma$ to $B$ such that 
for some permutation $\rho$ of $(1, 2, 3, 4, \dots, m)$ and some set of 
mutually different meta-name variables $\{u_1, u_2, \dots, u_s \}$  ($s \geq 0$, $n+s =m$) with 
$$\{x_1, x_2, x_3, x_4, \dots, x_n \} \cap \{u_1, u_2, \dots, u_s \} = \emptyset,$$
$$ \{y_1, y_2, y_3, y_4, \dots, y_m \} \cap \{u_1, u_2, \dots, u_s \} = \emptyset,$$
$$\sigma = \begin{pmatrix}
y_{\rho(1)} & y_{\rho(2)} & y_{\rho(3)} & y_{\rho(4)} & \dots & y_{\rho(n)} & y_{\rho(n+1)} & y_{\rho(n+2)} & \dots & y_{\rho(n+s)}  \\
x_1 & x_2 & x_3 & x_4 & \dots & x_n & u_1 & u_2 & \dots & u_s \\
\end{pmatrix}, $$

$$\sigma(B) \equiv A$$
is an instance of tautology of classical propositional logic.

(Case 2) When $n > m$, 
there is a uniform subsititution $\sigma$ to $A$ such that 
for some permutation $\rho$ of $(1, 2, 3, 4, \dots, n)$ and some set of 
mutually different meta-name variables $\{v_1, v_2, \dots, v_t \}$  ($t > 0$, $m+t =n$) with 
$$\{x_1, x_2, x_3, x_4, \dots, x_n \} \cap \{v_1, v_2, \dots, v_t \} = \emptyset,$$
$$ \{y_1, y_2, y_3, y_4, \dots, y_m \} \cap \{v_1, v_2, \dots, v_t \} = \emptyset,$$
$$\sigma = \begin{pmatrix}
x_{\rho(1)} & x_{\rho(2)} & x_{\rho(3)} & x_{\rho(4)} & \dots & x_{\rho(m)} & x_{\rho(m+1)} & x_{\rho(m+2)} & \dots & x_{\rho(m+t)}  \\
y_1 & y_2 & y_3 & y_4 & \dots & y_m & v_1 & v_2 & \dots & v_t \\
\end{pmatrix}, $$

$$\sigma(A) \equiv B$$
is an instance of tautology of classical propositional logic.

If $A$ is not quasi-trivial \rm (with respect to $B$\rm ), it is said to be \it quasi-nontrivial \rm (\it with respect to $B$\rm).   (The end of \bf Dfinition \ref{Definition 5}\rm )
\end{definition}
\medskip

Recall again  
\bigskip

($\mathrm{A_{M8}}$) $\enspace$ $\epsilon ab \wedge \epsilon cd . \supset . \epsilon aa \wedge \epsilon cc \wedge (\epsilon bc \supset . \epsilon ad \wedge \epsilon ba)$.

\bigskip

\begin{proposition} 
$\mathrm{(A_{M8})}$  is quasi-nontrivial \rm (\it with respect to $\mathrm{(A_t)}$\rm).
\end{proposition}
\medskip

Proof Trivial from Proposition \ref{Proposition 4} and Definition \ref{Definition 5}. $\Box$
\medskip

An example of quasi-trivial single axiom schemata with respect to $(\mathrm{A_{M8}})$ is, for example,

\medskip

$\enspace$ $\epsilon ab \wedge \epsilon de. \supset . \epsilon dd \wedge \epsilon aa \wedge (\epsilon bd \supset \epsilon ae) \wedge (\neg \epsilon ba \supset \neg \epsilon bd).  \enspace (*)$
\medskip

Take 
$$\sigma = \begin{pmatrix}
a & b & d & e \\
a & b & c & d \\
\end{pmatrix}. $$
Then we have $$\sigma((*)) \equiv (\mathrm{A_{M8}})$$ 
with $\rho = id$.

We shall give one more example of quasi-trivial single axiom schemata with respect to $(\mathrm{A_{M8}})$ is, for example,

\medskip

$\enspace$ $\epsilon ab \wedge \epsilon de. \supset . \epsilon dd \wedge \epsilon aa \wedge (\epsilon bd \supset . \epsilon ae \wedge \epsilon ba) \wedge 
(\epsilon cc \vee \neg \epsilon cc).  \enspace (**)$

\medskip

Take $v$ such that 
$$\{a, b, c, d, e \} \cap \{ v \} = \emptyset.$$
Take further 
$$\sigma = \begin{pmatrix}
a & b & d & e & c \\
a & b & c & d & v \\
\end{pmatrix}. $$

Then we have $$\sigma((**)) \equiv (\mathrm{A_{M8}})$$ 
with $\rho = id$.

We are not interested in axiom schemata characteristic of $\bf L_1$ which are quasi-trivial with each other. 
By the definition, we say that they are substantially the same ones.

\section{Structure theorems for nontriviality and quasi-nontriviality}
In this section, we shall discuss about the structures on nontriviality and quasi-nontriviality.
\medskip

\begin{definition}\label{Definition 6}
We define 

$(1)$ $Triv_{At}$  as a one-place predicate such that 
\smallskip

$\qquad \qquad Triv_{At}(x) \Leftrightarrow (x \enspace\mbox{is trivial with respect to} \enspace \mathrm{ (A_t)}),$

$(2)$ $NTriv_{At}$  as a one-place predicate such that 
\smallskip

$\qquad \qquad NTriv_{At}(x) \Leftrightarrow (x \enspace\mbox{is nontrivial with respect to} \enspace \mathrm{ (A_t)}).$
\end{definition}
\medskip

\begin{definition}\label{Definition 7}
We define 

$(1)$ $QTriv$  as a two-places predicate such that 
\smallskip

$\qquad \qquad QTriv(x, y) \Leftrightarrow (x \enspace\mbox{is quasi-trivial with respect to} \enspace y),$

$(2)$ $QNTriv$  as a two-places predicate such that 
\smallskip

$\qquad \qquad QNTriv(x, y) \Leftrightarrow (x \enspace\mbox{is quasi-nontrivial with respect to} \enspace y).$
\end{definition}
\medskip

\begin{definition}
We define 

$(1)$ $<_{QT}$  as a binary relation such that 
\smallskip

$\qquad \qquad  x <_{QT} y \Leftrightarrow QTriv(x, y),$

$(2)$ $<_{QNT}$  as a binary relation such that 
\smallskip

$\qquad \qquad  x <_{QNT} y \Leftrightarrow QNTriv(x, y).$
\end{definition}
\medskip

We note the following.
\medskip

\begin{proposition} We have

$(1)$ $Triv_{At}(\mathrm{ (A_t)}),$

$(2)$ $QTriv(x,  \mathrm{ (A_t)}) \Leftrightarrow Triv_{At}(x),$

$(3)$ $QNTriv(x,  \mathrm{ (A_t)}) \Leftrightarrow NTriv_{At}(x).$
\end{proposition}
\medskip

Proof. Trivial from Definitions \ref{Definition 5}, \ref{Definition 6} and \ref{Definition 7}.  $\Box$
\medskip

\begin{proposition} For any $x, y, z$ we have

$(1)$ $<_{QT}$ is reflexive and symmetric,

$(2)$ If $\# nv(x) \leq \#nv(y) \leq \#nv(z)$, $x <_{QT} y$ and $y <_{QT} z$, then we have $x <_{QT} z,$

$(3)$ If $\# nv(x) \geq \#nv(y) \geq \#nv(z)$, $x <_{QT} y$ and $y <_{QT} z$, then we have $x <_{QT} z$. 
\end{proposition}
\medskip

Proof. Easy from definitions. For (3), use (1). $\Box$
\medskip

That is , $<_{QT}$ is transitive for any monotonically increasing (decreasing) sequence $(x, y, z)$ with respect to $\#nv$.
\medskip

\begin{proposition}\label{Proposition 8} $<_{QNT}$ is symmetric.
\end{proposition}
\medskip

Proof. Easy from definitions. $\Box$

We note that $<_{QNT}$ is, in general, not transitive.

\section{A $\bf L_1$-provable formula to deduce transitivity (Ax2) and exchangeability (Ax3)}

We shall give a $\bf L_1$-provable formula to deduce transitivity (Ax2) and exchangeability (Ax3) as follows.

\bigskip

($\rm A_{S3}$) $\enspace$ $\epsilon ab \wedge \epsilon bc . \supset . \epsilon bb \wedge (\epsilon cd \supset . \epsilon ad \wedge \epsilon ba)$

\bigskip

\noindent We can easily show that ($\rm A_{S3}$) is provable in $\bf L_1$. From (Ax1) we have $$\epsilon ab \wedge \epsilon bc . \supset  \epsilon bb. \enspace \mbox{(5.1)}$$

\noindent From assumptions $\epsilon ab \wedge \epsilon bc$ and $\epsilon cd$ we obtain $\epsilon ad$ by applying (Ax2) twice and using propositional logic. So we have $$\epsilon ab \wedge \epsilon bc \wedge \epsilon cd. \supset \epsilon ad.\enspace \mbox{(5.2)}$$

\noindent By propositional logic and (Ax3),  $$\epsilon ab \wedge \epsilon bc \wedge \epsilon cd . \supset \epsilon ba \enspace \mbox{(5.3)}$$ holds. From (5.2) and (5.3), we get  $$\epsilon ab \wedge \epsilon bc \wedge \epsilon cd . \supset . \epsilon ad \wedge \epsilon ba \enspace \mbox{(5.4)}$$ Then we obtain ($\rm A_{S3}$) from (5.1) and (5.4) and propositional logic.

Let us see the following single axiom-schema characteristic of (Ax2) and (Ax3) of $\bf L_1$. 
\bigskip  

($\rm A_{t-1}$) $\enspace$ $\epsilon ab \wedge \epsilon bc . \supset . \epsilon ac \wedge \epsilon ba$.
\bigskip 

\noindent We easily see that $(A_{\rm t - 1}) \equiv . (Ax2) \wedge (Ax3)$ is an instance of tautology.
\bigskip

We can easily adopt the definition of non-triviality and quasi-nontriviality for ($\rm A_{S3}$) by replacing ($\rm A_{t}$) by ($\rm A_{t-1}$). We shall prove that  ($\rm A_{S3}$) is a non-trivial schema for (Ax2) and (Ax3) of $\bf L_1$ with respect to ($\rm A_{t-1}$).
\medskip 

\begin{proposition} 
$\mathrm{(A_{S3})}$  is nontrivial \rm (\it with respect to $\mathrm{(A_{t - 1}))}$.
\end{proposition}
\medskip

Proof.
We shall show $\mathrm{(A_{S3})}$ is nontrivial. We see 
$$nv(\mathrm{A_{S3}}) = (a, b, c, d), \enspace \# nv(\mathrm{A_{S3}})= 4.$$ 

The idea is to give a sentential valuation $v$ such that the premise of $\sigma(\rm A_{S3})$ has the value $t$ and its conclusion has $f$ with respect to 
each permutation $\sigma$ of meta-name variables, whereas $v(\mathrm{(A_{t-1})}) = t$ holds under the $v$.

Take a meta-name variable $y$ such that $\{a, b, c, d\} \cap \{y\} = \emptyset$.
\medskip

(Case 1) Let
$$\sigma = \begin{pmatrix}
u & b & v & w \\
a & b & c & y \\
\end{pmatrix}.$$

Then take a sentential valuation $v$ such that $v(\epsilon bb) = f$. For the rest $t$ is assigned.

Then we easily see $v((\rm A_{t-1})) = \it t$ and the premise of $v(\sigma(\rm A_{S3})) = \it f$ under the assignment $v$. Thus $\sigma(\rm A_{S3}) \equiv (\rm A_{t-1})$ is not an instance of tautology. 
\medskip


(Case 2) Let $x \neq b$. 

\noindent (Subcase 2.1) Let 
$$\sigma = \begin{pmatrix}
b & x& v & w \\
a & b & c & y \\
\end{pmatrix}.$$

Then take a sentential valuation $v$ such that $v(\epsilon aa) = f$. For the rest $t$ is assigned.

Then we easily see $v((\rm A_{t-1})) = \it t$ and the premise of $v(\sigma(\rm A_{S3})) = \it f$ under the assignment $v$. So $\sigma(\rm A_{S3}) \equiv (\rm A_{t-1})$ is not an instance of tautology. 



\noindent (Subcase 2.2) Let 
$$\sigma = \begin{pmatrix}
u & x & b & w \\

a & b & c & y \\
\end{pmatrix}.$$

Then take a sentential valuation $v$ such that $v(\epsilon cc) = f$. For the rest $t$ is assigned.

Then we easily see $v((\rm A_{t-1})) = \it t$ and the premise of $v(\sigma(\rm A_{S3})) = \it f$ under the assignment $v$. Hence $\sigma(\rm A_{S3}) \equiv (\rm A_{t-1})$ is not an instance of tautology. 



\noindent (Subcase 2.3) Let 
$$\sigma = \begin{pmatrix}
u & x & v & b \\
a & b & c & y \\
\end{pmatrix}.$$

Then take a sentential valuation $v$ such that $v(\epsilon yy) = f$. For the rest is assigned to $t$.

Then we easily see $v((\rm A_{t-1})) = \it t$ and the premise of $v(\sigma(\rm A_{S3})) = \it f$ under the assignment $v$, Thus $\sigma(\rm A_{S3}) \equiv (\rm A_{t-1})$ is not an instance of tautology.  $\Box$

Lastly, we shall show that (Ax2) and (Ax3) are deducible from ($\rm A_{S3}$). 
\medskip

\begin{proposition} 
$(\rm Ax2)$ and $(\rm Ax3)$ are deducible from $(\rm A_{S3})$. 
\end{proposition}
\medskip

\noindent Proof.

Assume 

($\rm A_{S3}$) $\enspace$ $\epsilon ab \wedge \epsilon bc . \supset . \epsilon bb \wedge (\epsilon cd \supset . \epsilon ad \wedge \epsilon ba).$
\medskip

\noindent From ($\rm A_{S3}$) we have
$$\epsilon ab \wedge \epsilon bc . \supset \epsilon bb. \quad \mbox{(5.4)}$$
$$\epsilon ab \wedge \epsilon bc . \supset (  \epsilon cd \supset . \epsilon ad \wedge \epsilon ba). \quad \mbox{(5.5)}$$
By propositional logic and (5.5) we obtain 
$$\epsilon ab \wedge \epsilon cd . \supset (  \epsilon bc \supset . \epsilon ad \wedge \epsilon ba). \quad \mbox{(5.6)}$$
By substitution with $b \rightarrow c, c \rightarrow d$ in (5.6) and Proposition \ref{Proposition 1}, we get 
$$\epsilon ab \wedge \epsilon bc . \supset (  \epsilon bb \supset . \epsilon ac \wedge \epsilon ba). \quad \mbox{(5.7)}$$ 
Here, set 
$$A =_{def.} \epsilon ab \wedge \epsilon bc,$$
$$B =_{def.} \epsilon bb,$$
$$C =_{def.} \epsilon ac \wedge \epsilon ba.$$
We know that 
$$(A \supset B) \wedge (A \supset . B \supset C ). \supset . A \supset C. \quad \mbox{(5.8)}$$ 
is a tautology of classical propositional logic. So from (5.4), (5,7) and (5,8), we get $A \supset C$, that is ($\rm A_{t-1}$). Thus (Ax2) and (Ax3) hold. $\Box$

\section{Simplified axiom schemata $\mathrm{(A_{S1})}$, $\mathrm{(A_{S2})}$, $\mathrm{(A_{S3N})}$ and $\mathrm{(A_{S3Nd})}$ characteristic of $\bf L_1$}

We shall now introduce three single axiom schemata characteristic of $\bf L_1$ as simplified ones of ($\rm A_{M8}$) as follows:

\bigskip 

($\rm A_{S1}$) $\enspace$ $\epsilon ab \wedge \epsilon cd . \supset . \epsilon aa \wedge (\epsilon bc \supset . \epsilon ad \wedge \epsilon ba)$.

\bigskip 

($\rm A_{S2}$) $\enspace$ $\epsilon ab \wedge \epsilon cd . \supset . \epsilon cc \wedge (\epsilon bc \supset . \epsilon ad \wedge \epsilon ba)$.

\bigskip 

($\rm A_{S3N}$) $\enspace$ $\epsilon ab \supset .\epsilon aa \wedge (\epsilon bc \supset . \epsilon bb \wedge (\epsilon cd \supset . \epsilon ad \wedge \epsilon ba))$. 

\bigskip 

($\rm A_{S3Nd}$) $\enspace$ $\epsilon ab \supset .\epsilon aa \wedge (\epsilon bc \wedge \epsilon cd . \supset . \epsilon ad \wedge \epsilon ba)$. 
\medskip

\begin{proposition} $ $
\medskip

$(1)$  $\mathrm{(A_{S1})}$ $\Longleftrightarrow$ $\bf L_1$ 
\medskip

$(2)$  $\mathrm{(A_{S2})}$ $\Longleftrightarrow$ $\bf L_1$
\medskip

$(3)$  $\mathrm{(A_{S3N})}$ $\Longleftrightarrow$ $\bf L_1$
\medskip

$(4)$  $\mathrm{(A_{S3Nd})}$ $\Longleftrightarrow$ $\bf L_1$
\medskip
\end{proposition}

\noindent Proof
\bigskip

\noindent The proof of (1). 

(The case of $\Longrightarrow$) 
Take the following uniform substitutions:

$$\sigma_1 = \begin{pmatrix}
a & b & c & d \\
a & b & a & b \\
\end{pmatrix},$$

$$\sigma_2 = \begin{pmatrix}
a & b & c & d \\
a & b & b & c \\
\end{pmatrix}.$$

By propositional logic, we immediately derive (Ax1) from $\sigma_1(\rm A_{S1})$. By the obtained (Ax1) and propositional logic, we immediately derive (Ax2) and (Ax3) from $\sigma_2(\rm A_{S1})$. 
\medskip

(The case of $\Longleftarrow$) $(\rm A_{S1})$ is a theorem of $\bf L_1$. We can directly derive ($\rm A_{S1}$) from (Ax1)--(Ax3). By (Ax1), we have 
\medskip

(B1) $\epsilon ab \supset \epsilon aa$.
\medskip

\noindent From (B1) we have
\medskip

(B2) $\epsilon ab \wedge \epsilon cd . \supset \epsilon aa$.

\medskip

\noindent  By using (Ax2) two times and (Ax3) and propositional logic, we easily get

\medskip

(B3) $\epsilon ab \wedge \epsilon bc \wedge \epsilon cd. \supset \epsilon ad$, 

\medskip

(B4) $\epsilon ab \wedge \epsilon bc \wedge \epsilon cd. \supset \epsilon ba$.

\medskip

\noindent Then we can derive $ (\rm A_{S1})$ from (B2), (B3) and (B4).
\medskip

\noindent The proof of (2).

(The case of $\Longrightarrow$) 
Take the following uniform substitutions:

$$\sigma_1 = \begin{pmatrix}
a & b & c & d \\
a & b & a & b \\
\end{pmatrix},$$

$$\sigma_2 = \begin{pmatrix}
a & b & c & d \\
a & b & b & c \\
\end{pmatrix}.$$

By propositional logic, we immediately derive (Ax1) from $\sigma_1(\rm A_{S2})$. By the obtained (Ax1) and propositional logic, we immediately derive (Ax2) and (Ax3) from $\sigma_2(\rm A_{S2})$. 
\medskip

(The case of $\Longleftarrow$) $(\rm A_{S2})$ is a theorem of $\bf L_1$. We can directly derive ($\rm A_{S2}$) from (Ax1)--(Ax3). By (Ax1), we have 
\medskip

(C1) $\epsilon cd \supset \epsilon cc$.
\medskip

\noindent From (C1) we have
\medskip

(C2) $\epsilon ab \wedge \epsilon cd . \supset \epsilon cc$.

\medskip

\noindent  By using (Ax2) two times and (Ax3) and propositional logic, we easily get

\medskip

(C3) $\epsilon ab \wedge \epsilon bc \wedge \epsilon cd. \supset \epsilon ad$, 

\medskip

(C4) $\epsilon ab \wedge \epsilon bc \wedge \epsilon cd. \supset \epsilon ba$.

\medskip

\noindent Then we can derive $ (\rm A_{S2})$ from (C2), (C3) and (C4).
\medskip

\noindent The proof of (3).

(The case of $\Longrightarrow$) 
Take the following uniform substitutions:

$$\sigma_1 = \begin{pmatrix}
a & b & c & d \\
a & b & a & b \\
\end{pmatrix},$$

$$\sigma_2 = \begin{pmatrix}
a & b & c & d \\
a & b & b & c \\
\end{pmatrix}.$$

By propositional logic, we immediately derive (Ax1) from $\sigma_1(\rm A_{S3N})$. By the obtained (Ax1) and propositional logic, we immediately derive (Ax2) and (Ax3) from $\sigma_2(\rm A_{S3N})$.
\medskip 

(The case of $\Longleftarrow$) $(\rm A_{S3N})$ is a theorem of $\bf L_1$. We can directly derive ($\rm A_{S3N}$) from (Ax1)--(Ax3). By (Ax1), we have 
\medskip

(D1a) $\epsilon ab \supset \epsilon aa$, 
\medskip

(D1b) $\epsilon bc \supset \epsilon bb$.
\medskip

\noindent By using (Ax2) two times and (Ax3) and propositional logic, we easily get

\medskip

(D2) $\epsilon ab \wedge \epsilon bc \wedge \epsilon cd. \supset \epsilon ad$, 

\medskip

(D3) $\epsilon ab \wedge \epsilon bc \wedge \epsilon cd. \supset \epsilon ba$.

\medskip

\noindent Then we can derive $ (\rm A_{S3N})$ from (D1a), (D1b), (D2) and (D3). 

\medskip

\noindent The proof of (4). 

(The case of $\Longrightarrow$) 
Take the following uniform substitutions:

$$\sigma_1 = \begin{pmatrix}
a & b & c & d \\
a & b & a & b \\
\end{pmatrix},$$

$$\sigma_2 = \begin{pmatrix}
a & b & c & d \\
a & b & b & c \\
\end{pmatrix}.$$

Then $\epsilon ab \supset \epsilon aa. \equiv \sigma_1(\rm A_{S3Nd})$ is a tautology of propositional logic. So from this, we have (Ax1).

By the obtained (Ax1) and propositional logic, we immediately get 
\medskip

(E1) $\epsilon bc \supset \epsilon bb.$
\medskip

\noindent By $\sigma_2(\rm A_{S3Nd})$ and propositional logic, we have 
\medskip

(E2) $\epsilon ab \supset .  \epsilon bb \wedge \epsilon bc. \supset . \epsilon ac \wedge  \epsilon ba$. 
\medskip

\noindent So by (E2) and propositional logic we get
\medskip
 
(E3) $\epsilon ab \wedge \epsilon bb \wedge \epsilon bc. \supset . \epsilon ac \wedge  \epsilon ba$.
\medskip

\noindent By (E1) and proposisional logic, we obtain 
\medskip

(E4) $\epsilon ab \wedge \epsilon bb \wedge \epsilon bc. \supset . \epsilon ac \wedge  \epsilon ba$.
\medskip

\noindent By (E1) and (E3) and proposisional logic, we obtain 
\medskip
 
(E5) $\epsilon ab \wedge \epsilon bb \wedge \epsilon bc. \supset . \epsilon ac \wedge  \epsilon ba$.
\medskip

derive (Ax2) and (Ax3) from $\sigma_2(\rm A_{S3Nd})$.
\medskip 

(The case of $\Longleftarrow$) $(\rm A_{S3Nd})$ is a theorem of $\bf L_1$. We can directly derive ($\rm A_{S3Nd}$) from (Ax1)--(Ax3). By (Ax1), we have 
\medskip

(F1) $\epsilon ab \supset \epsilon aa$, 
\medskip

\noindent By using (Ax2) two times and (Ax3) and propositional logic, we easily get

\medskip

(F2) $\epsilon ab \wedge \epsilon bc \wedge \epsilon cd. \supset \epsilon ad$, 

\medskip

(F3) $\epsilon ab \wedge \epsilon bc \wedge \epsilon cd. \supset \epsilon ba$.

\medskip

\noindent Then we can derive $ (\rm A_{S3Nd})$ from (F1), (F2) and (F3). $\enspace \Box$

\section{Nontriviality of $\mathrm{(A_{S1})}$, $\mathrm{(A_{S2})}$, $\mathrm{(A_{S3N})}$ and $\mathrm{(A_{S3Nd})}$ with respect to $\mathrm{(A_{t})}$} 

\begin{proposition}\label{Proposition 12} We have:

$(1)$ $\mathrm{(A_{S1})}$  is nontrivial with respect to $\mathrm{(A_t)}$.

$(2)$ $\mathrm{(A_{S2})}$  is nontrivial with respect to $\mathrm{(A_t)}$.

$(3)$ $\mathrm{(A_{S3N})}$  is nontrivial with respect to $\mathrm{(A_t)}$.

$(4)$ $\mathrm{(A_{S3Nd})}$  is nontrivial with respect to $\mathrm{(A_t)}$.

\end{proposition}
\medskip

\noindent Proof. One can easily see that ($\rm A_{S1}$), ($\rm A_{S2}$) and $\mathrm{(A_{S3N})}$ are nontrivial single axiom schemata characteristic of $\bf L_1$, respectively, as in the proof of Proposition \ref{Proposition 4}.

\noindent The proof of (1).

We shall show that is nontrivial. We see 
$$nv(\mathrm{A_{S1}}) = (a, b, c, d), \enspace \# nv(\mathrm{A_{S1}})= 4.$$ 

In order to prove the nontriviality of ($\rm A_{S1}$), we must check 24 ($= 4!$) cases, that is, 
the number of permutations of $(1, 2, 3, 4)$. However, the real case to be checked is a few.

Take a meta-name variable $y$ such that $\{a, b, c, d\} \cap \{y\} = \emptyset$. We use $t$ (= true) and $f$ (= false) as usual.

\medskip

(Case 1) Let
$$\sigma = \begin{pmatrix}
a & u & v & w \\
a & b & c & y \\
\end{pmatrix}.$$

In this case, the idea is to give a sentential valuation $v$ such that, noticing the form as $\epsilon \alpha \alpha$, the premise of $\mathrm{(A_t)}$ has the value $t$ and the conclution of $\mathrm{(A_t)}$ has the value $f$ with respect to 
each permutation $\sigma$ of meta-name variables, whereas the premise of $\mathrm{(A_{S1})}$ has the valure $f$ under the $v$.

Then take a sentential valuation $v$ such that $v(\epsilon aa) = v(\epsilon ac) = v(\epsilon ay) = f$. For the rest $t$ is assigned.

Then we easily see $v((\rm A_t)) = \it f$ and $v(\sigma(\rm A_{S1})) = \it t$ under the assignment $v$. Thus $\sigma(\rm A_{S1}) \equiv (\rm A_t)$ is not an instance of tautology. 
\medskip


(Case 2) Let $x \neq a$. 

\noindent (Subcase 2.1) Let 
$$\sigma = \begin{pmatrix}
x & a& v & w \\
a & b & c & y \\
\end{pmatrix}.$$

In this case, the idea is to give a sentential valuation $v$ such that, noticing the form as $\epsilon \alpha \alpha$, the premise of $\mathrm{(A_{t})}$ has the value $t$ and the conclution of $\mathrm{(A_{t})}$ has the value $f$ with respect to 
each permutation $\sigma$ of meta-name variables, whereas the premise of $\mathrm{(A_{S1})}$ has the valure $f$ under the $v$.

Then take a sentential valuation $v$ such that $v(\epsilon bb) = f$. For the rest $t$ is assigned.

Then we easily see $v((\rm A_t)) = \it t$ and $v(\sigma(\rm A_{S1})) = \it f$ under the assignment $v$. Thus $\sigma(\rm A_{S1}) \equiv (\rm A_t)$ is not an instance of tautology. 
\medskip



\noindent (Subcase 2.2) Let 
$$\sigma = \begin{pmatrix}
x & u & a & w \\
a & b & c & y \\
\end{pmatrix}.$$

Then take a sentential valuation $v$ such that $v(\epsilon cc) = f$. For the rest $t$ is assigned.

Then we easily see $v((\rm A_t)) = \it t$ and $v(\sigma(\rm A_{S1})) = \it f$ under the assignment $v$. Hence $\sigma(\rm A_{S1}) \equiv (\rm A_t)$ is not an instance of tautology. 



\noindent (Subcase 2.3) Let 
$$\sigma = \begin{pmatrix}
x & u & v & a \\
a & b & c & y \\
\end{pmatrix}.$$

Then take a sentential valuation $v$ such that $v(\epsilon yy) = f$. For the rest is assigned to $t$.

Then we easily see $v((\rm A_t)) = \it t$ and $v(\sigma(\rm A_{S1})) = \it f$ under the assignment $v$, Thus $\sigma(\rm A_{S1}) \equiv (\rm A_t)$ is not an instance of tautology.  $\Box$
\medskip

\noindent The proof of (2).

We shall show that is nontrivial. We see 
$$nv(\mathrm{A_{S2}}) = (a, b, c, d), \enspace \# nv(\mathrm{A_{S2}})= 4.$$ 
. 
%
(Case 1) Let
$$\sigma = \begin{pmatrix}
u & v & c & w \\
a & b & c & y \\
\end{pmatrix}.$$

Then take a sentential valuation $v$ such that $v(\epsilon cc) = f$. For the rest $t$ is assigned.

Then we easily see $v((\rm A_t)) = \it t$ and $v(\sigma(\rm A_{S2})) = \it f$ under the assignment $v$. Thus $\sigma(\rm A_{S2}) \equiv (\rm A_t)$ is not an instance of tautology. 
\medskip


(Case 2) Let $x \neq c$. 

\noindent (Subcase 2.1) Let 
$$\sigma = \begin{pmatrix}
c & v & x & w \\
a & b & c & y \\
\end{pmatrix}.$$

Then take a sentential valuation $v$ such that $v(\epsilon aa) = v(\epsilon ac) = v(\epsilon ay) = v(\epsilon cy) = v(\epsilon yc) = f$ and $v(\epsilon ab) = t$. For the rest $t$ is assigned.

Then we easily see $v((\rm A_t)) = \it f$ and $v(\sigma(\rm A_{S2})) = \it t$ under the assignment $v$. So $\sigma(\rm A_{S2}) \equiv (\rm A_t)$ is not an instance of tautology. 



\noindent (Subcase 2.2) Let 
$$\sigma = \begin{pmatrix}
u & c & x & w \\
a & b & c & y \\
\end{pmatrix}.$$

Then take a sentential valuation $v$ such that $v(\epsilon bb) = f$. For the rest $t$ is assigned.

Then we easily see $v((\rm A_t)) = \it t$ and $v(\sigma(\rm A_{S2})) = \it f$ under the assignment $v$. Hence $\sigma(\rm A_{M8}) \equiv (\rm A_t)$ is not an instance of tautology. 



\noindent (Subcase 2.3) Let 
$$\sigma = \begin{pmatrix}
u & v & x & c \\
a & b & c & y \\
\end{pmatrix}.$$

Then take a sentential valuation $v$ such that $v(\epsilon yy) = f$. For the rest is assigned to $t$.

Then we easily see $v((\rm A_t)) = \it t$ and $v(\sigma(\rm A_{S2})) = \it f$ under the assignment $v$, Thus $\sigma(\rm A_{M8}) \equiv (\rm A_t)$ is not an instance of tautology.  $\Box$
\medskip

\noindent The proof of (3).

We see 
$$nv(\mathrm{A_{S3N}}) = (a, b, c, d), \enspace \# nv(\mathrm{A_{S3N}})= 4.$$ 

The idea is to give a sentential valuation $v$ such that the premise of $\sigma(\rm A_{S3N})$ has the value $t$ and its conclusion has $f$ with respect to 
each permutation $\sigma$ of meta-name variables, whereas $v(\mathrm{(A_t)}) = t$ holds under the $v$.

Take a meta-name variable $y$ such that $\{a, b, c, d\} \cap \{y\} = \emptyset$. We use $t$ (= true) and $f$ (= false) as usual.
\medskip

(Case 1) Let
$$\sigma = \begin{pmatrix}
u & b & v & w \\
a & b & c & y \\
\end{pmatrix}.$$

Then take a sentential valuation $v$ such that $v(\epsilon bb) = f$. For the rest $t$ is assigned.

Then we easily see $v((\rm A_t)) = \it t$ and $v(\sigma(\rm A_{S3N})) = \it f$ under the assignment $v$. Thus $\sigma(\rm A_{S3N}) \equiv (\rm A_t)$ is not an instance of tautology. 
\medskip


(Case 2) Let $x \neq b$. 

\noindent (Subcase 2.1) Let 
$$\sigma = \begin{pmatrix}
b & x & v & w \\
a & b & c & y \\
\end{pmatrix}.$$

Then take a sentential valuation $v$ such that $v(\epsilon bb) = v(\epsilon cc) = v(\epsilon yy) = f$. For the rest $t$ is assigned.

Then we easily see $v((\rm A_t)) = \it t$ and $v(\sigma(\rm A_{S3N})) = \it f$ under the assignment $v$. So $\sigma(\rm A_{S3N}) \equiv (\rm A_t)$ is not an instance of tautology. 



\noindent (Subcase 2.2) Let 
$$\sigma = \begin{pmatrix}
u & x & b & w \\
a & b & c & y \\
\end{pmatrix}.$$

Then take a sentential valuation $v$ such that $v(\epsilon cc) = f$. For the rest $t$ is assigned.

Then we easily see $v((\rm A_t)) = \it t$ and $v(\sigma(\rm A_{S3N})) = \it f$ under the assignment $v$. Hence $\sigma(\rm A_{S3N}) \equiv (\rm A_t)$ is not an instance of tautology. 



\noindent (Subcase 2.3) Let 
$$\sigma = \begin{pmatrix}
u & x & v & b \\
a & b & c & y \\
\end{pmatrix}.$$

Then take a sentential valuation $v$ such that $v(\epsilon yy) = f$. For the rest is assigned to $t$.

Then we easily see $v((\rm A_t)) = \it t$ and $v(\sigma(\rm A_{S3N})) = \it f$ under the assignment $v$, Thus $\sigma(\rm A_{S3N}) \equiv (\rm A_t)$ is not an instance of tautology. 
\medskip

\noindent The proof of (4).

We see 
$$nv(\mathrm{A_{S3Nd}}) = (a, b, c, d), \enspace \# nv(\mathrm{A_{S3Nd}})= 4.$$ 

The idea is to give a sentential valuation $v$ such that the premise of $\sigma(\rm A_{S3Nd})$ has the value $t$ and its conclusion has $f$ with respect to 
each permutation $\sigma$ of meta-name variables, whereas $v(\mathrm{(A_t)}) = t$ holds under the $v$.

Take a meta-name variable $y$ such that $\{a, b, c, d\} \cap \{y\} = \emptyset$. 
\medskip

(Case 1) Let
$$\sigma = \begin{pmatrix}
a & b & c & d \\
a & b & c & y \\
\end{pmatrix}.$$

\noindent Then take a sentential valuation $v$ such that $v(\epsilon ay) = f$. For the rest $t$ is assigned.

Then we easily see $v((\rm A_t)) = \it t$ and $v(\sigma(\rm A_{S3Nd})) = \it f$ under the assignment $v$. Thus $\sigma(\rm A_{S3Nd}) \equiv (\rm A_t)$ is not an instance of tautology. 
\medskip


(Case 2) Let 
$$\sigma = \begin{pmatrix}
a & b & d & c \\
a & b & c & y \\
\end{pmatrix}.$$

\noindent Then take a sentential valuation $v$ such that $v(\epsilon bc) = v(\epsilon ba) = f$. For the rest $t$ is assigned.

Then we easily see $v((\rm A_t)) = \it t$ and $v(\sigma(\rm A_{S3Nd})) = \it f$ under the assignment $v$. So $\sigma(\rm A_{S3Nd}) \equiv (\rm A_t)$ is not an instance of tautology. 



\noindent (Case 3) Let $u \neq b$.

Let 
$$\sigma = \begin{pmatrix}
u & x & b & w \\
a & b & c & y \\
\end{pmatrix}.$$

\noindent Then take a sentential valuation $v$ such that $v(\epsilon ca) = v(\epsilon ya) =f$. For the rest $t$ is assigned.

Then we easily see $v((\rm A_t)) = \it t$ and $v(\sigma(\rm A_{S3Nd})) = \it f$ under the assignment $v$. Hence $\sigma(\rm A_{S3Nd}) \equiv (\rm A_t)$ is not an instance of tautology. 



\noindent (Case 4)  Let $u \neq a$

Let 
$$\sigma = \begin{pmatrix}
u & x & v & b \\
a & b & c & y \\
\end{pmatrix}.$$

\noindent Then take a sentential valuation $v$ such that $v(\epsilon bb) = v(\epsilon cc) = v(\epsilon yy) = f$. For the rest is assigned to $t$.

Then we easily see $v((\rm A_t)) = \it t$ and $v(\sigma(\rm A_{S3Nd})) = \it f$ under the assignment $v$, Thus $\sigma(\rm A_{S3Nd}) \equiv (\rm A_t)$ is not an instance of tautology.  $\Box$

\section{The quasi-nontrivialities}

We shall proceed to see natural relationships among the proposed axiom schemata by means of the quasi-nontriviality.
\medskip 

\begin{theorem} Four axiom schemata $(\rm A_{M8})$, $(\rm A_{S1})$, $(\rm A_{S2})$, $(\rm A_{S3N})$ and $(\rm A_{S3Nd})$ characteristic of $\bf L_1$ 
are nontrivial and further we know the following:
	
	$(1)$ $(\rm A_{M8})$ is quasi-nontrivial with respect to  $(\rm A_{S1})$.
	
	$(2)$ $(\rm A_{M8})$ is quasi-nontrivial with respect to  $(\rm A_{S2})$.
	
	$(3)$ $(\rm A_{M8})$ is quasi-nontrivial with respect to  $(\rm A_{S3N})$.
	
	$\rm (3d)$ $(\rm A_{M8})$ is quasi-nontrivial with respect to  $(\rm A_{S3Nd})$.
	
	$(4)$ $(\rm A_{S1})$ is quasi-nontrivial with respect to  $(\rm A_{M8})$.
	
	$(5)$ $(\rm A_{S1})$ is quasi-nontrivial with respect to  $(\rm A_{S2})$.
	
	$(6)$ $(\rm A_{S1})$ is quasi-nontrivial with respect to  $(\rm A_{S3N})$
	
	$\rm (6d)$ $(\rm A_{S1})$ is quasi-nontrivial with respect to  $(\rm A_{S3Nd})$.
	
	$(7)$ $(\rm A_{S2})$ is quasi-nontrivial with respect to  $(\rm A_{M8})$.
	
	$(8)$ $(\rm A_{S2})$ is quasi-nontrivial with respect to  $(\rm A_{S1})$.
	
	$(9)$ $(\rm A_{S2})$ is quasi-nontrivial with respect to  $(\rm A_{S3N})$.
	
	$\rm (9d)$ $(\rm A_{S2})$ is quasi-nontrivial with respect to  $(\rm A_{S3Nd})$.
	
	$(10)$ $(\rm A_{S3N})$ is quasi-nontrivial with respect to  $(\rm A_{M8})$.
	
	$(11)$ $(\rm A_{S3N})$ is quasi-nontrivial with respect to  $(\rm A_{S1})$.
	
	$(12)$ $(\rm A_{S3N})$ is quasi-nontrivial with respect to  $(\rm A_{S2})$.
	
	$\rm (12d)$ $(\rm A_{S3N})$ is quasi-nontrivial with respect to $(\rm A_{S3Nd})$.
	
	$(13)$ $(\rm A_{S3Nd})$ is quasi-nontrivial with respect to  $(\rm A_{M8})$.
	
	$(14)$ $(\rm A_{S3Nd})$ is quasi-nontrivial with respect to  $(\rm A_{S1})$.
	
	$(15)$ $(\rm A_{S3Nd})$ is quasi-nontrivial with respect to  $(\rm A_{S2})$.
	
	$(16)$ $(\rm A_{S3Nd})$ is quasi-nontrivial with respect to $(\rm A_{S3N})$.
	
\end{theorem}
\medskip

Proof. Because of Propositon \ref{Proposition 8}, we may prove (1), (2), (3), (5), (6), (9), (3d), (6d), (9d) and (12d). Recall also the nontriviality of 
$(\rm A_{M8})$, $(\rm A_{S1})$, $(\rm A_{S2})$, $(\rm A_{S3N})$ and $(\rm A_{S3Nd})$ from Propositions \ref{Proposition 4} and \ref{Proposition 12}.

We shall first prove (1).  Recall 

\bigskip

($\mathrm{A_{M8}}$) $\enspace$ $\epsilon ab \wedge \epsilon cd . \supset . \epsilon aa \wedge \epsilon cc \wedge (\epsilon bc \supset . \epsilon ad \wedge \epsilon ba)$,

\bigskip

($\rm A_{S1}$) $\enspace$ $\epsilon ab \wedge \epsilon cd . \supset . \epsilon aa \wedge (\epsilon bc \supset . \epsilon ad \wedge \epsilon ba)$.

\bigskip

We have two cases for the proof 

(Case 1)  Let
$$\sigma = \begin{pmatrix}
x & y & z & w \\
a & b & c & d \\
\end{pmatrix}$$
with $z \neq a$. Take a sentential valuation $v$ such that $v(\epsilon cc) =  f$ and the rest is $t$. Then $v((\rm A_{M8}))$ = $f$ and $v(\sigma((\rm A_{S1})))$ = $t$.  

(Case 2) Let 
$$\sigma = \begin{pmatrix}
x & y & a & w \\
a & b & c & d \\
\end{pmatrix}.$$
Take a sentential valuation $v$ such that $v(\epsilon aa) =  f$ and the rest is $t$.  Then $v((\rm A_{M8}))$ = $f$ and $v(\sigma((\rm A_{S1})))$ = $t$.  

We shall prove (2). Recall 

\bigskip

($\mathrm{A_{M8}}$) $\enspace$ $\epsilon ab \wedge \epsilon cd . \supset . \epsilon aa \wedge \epsilon cc \wedge (\epsilon bc \supset . \epsilon ad \wedge \epsilon ba)$,

\bigskip 

($\rm A_{S2}$) $\enspace$ $\epsilon ab \wedge \epsilon cd . \supset . \epsilon cc \wedge (\epsilon bc \supset . \epsilon ad \wedge \epsilon ba)$.

\bigskip 

We have two cases for the proof 

(Case 1)  Let
$$\sigma = \begin{pmatrix}
x & y & z & w \\
a & b & c & d \\
\end{pmatrix}$$
with $x \neq c$. Take a sentential valuation $v$ such that $v(\epsilon aa) =  f$ and the rest is $t$. Then $v((\rm A_{M8}))$ = $f$ and $v(\sigma((\rm A_{S2})))$ = $t$.  

(Case 2) Let 
$$\sigma = \begin{pmatrix}
c& y & z & w \\
a & b & c & d \\
\end{pmatrix}.$$
Take a sentential valuation $v$ such that $v(\epsilon cc) =  f$ and the rest is $t$.  Then $v((\rm A_{M8}))$ = $f$ and $v(\sigma((\rm A_{S2})))$ = $t$. 

We shall prove (3). Recall 

\bigskip

($\mathrm{A_{M8}}$) $\enspace$ $\epsilon ab \wedge \epsilon cd . \supset . \epsilon aa \wedge \epsilon cc \wedge (\epsilon bc \supset . \epsilon ad \wedge \epsilon ba)$,

\bigskip 

($\rm A_{S3N}$) $\enspace$ $\epsilon ab \supset .\epsilon aa \wedge (\epsilon bc \supset . \epsilon bb \wedge (\epsilon cd \supset . \epsilon ad \wedge \epsilon ba))$. 

\bigskip

We also have two cases for that. 

(Case 1)  Let
$$\sigma = \begin{pmatrix}
x & y & z & w \\
a & b & c & d \\
\end{pmatrix}$$
\noindent with $y \neq c$.
 
\noindent (Subcase 1.1)  Let
$$\sigma = \begin{pmatrix}
x & a & z & w \\
a & b & c & d \\
\end{pmatrix}.$$
Take a sentential valuation $v$ such that $v(\epsilon bb) =  f$ and the rest is $t$. Then $v((\rm A_{M8}))$ = $t$ and $v(\sigma((\rm A_{S3N})))$ = $f$.

\noindent (Subcase 1.2)  Let
$$\sigma = \begin{pmatrix}
x & b & z & w \\
a & b & c & d \\
\end{pmatrix}.$$
Take a sentential valuation $v$ such that $v(\epsilon bb) =  f$ and the rest is $t$. Then $v((\rm A_{M8}))$ = $t$ and $v(\sigma((\rm A_{S3N})))$ = $f$.

\noindent (Subcase 1.3)  Let
$$\sigma = \begin{pmatrix}
x & d & z & w \\
a & b & c & d \\
\end{pmatrix}.$$
Take a sentential valuation $v$ such that $v(\epsilon dd) = v(\epsilon ca) = v(\epsilon cb) =  f$ and the rest is $t$. Then $v((\rm A_{M8}))$ = $t$ and $v(\sigma((\rm A_{S3N})))$ = $f$.

(Case 2) Let 
$$\sigma = \begin{pmatrix}
x & c & z & w \\
a & b & c & d \\
\end{pmatrix}.$$
\noindent (Subcase 2.1) Let 
$$\sigma = \begin{pmatrix}
b & c & z & w \\
a & b & c & d \\
\end{pmatrix}.$$
Take a sentential valuation $v$ such that $v(\epsilon aa) = v(\epsilon ba) = v(\epsilon ca) = v(\epsilon da) =  f$ and the rest is $t$.  Then $v((\rm A_{M8}))$ = $f$ and $v(\sigma((\rm A_{S3N})))$ = $t$. 

\noindent (Subcase 2.2) Let
$$\sigma = \begin{pmatrix}
x & c & b & w \\
a & b & c & d \\
\end{pmatrix}.$$
Take a sentential valuation $v$ such that $v(\epsilon cc) = v(\epsilon bc) = v(\epsilon ac) = v(\epsilon dc) =  f$ and the rest is $t$.  Then $v((\rm A_{M8}))$ = $f$ and $v(\sigma((\rm A_{S3N})))$ = $t$. 

\noindent (Subcase 2.3) Let
$$\sigma = \begin{pmatrix}
x & c & z & b \\
a & b & c & d \\
\end{pmatrix}.$$
Take a sentential valuation $v$ such that $v(\epsilon dd) = f$ and the rest is $t$.  Then $v((\rm A_{M8}))$ = $t$ and $v(\sigma((\rm A_{S3N})))$ = $f$. 

We shall prove (5). Recall 

\bigskip 

($\rm A_{S1}$) $\enspace$ $\epsilon ab \wedge \epsilon cd . \supset . \epsilon aa \wedge (\epsilon bc \supset . \epsilon ad \wedge \epsilon ba)$,

\bigskip 

($\rm A_{S2}$) $\enspace$ $\epsilon ab \wedge \epsilon cd . \supset . \epsilon cc \wedge (\epsilon bc \supset . \epsilon ad \wedge \epsilon ba)$.

\bigskip

We also have four cases for that. 

(Case 1)  Let
$$\sigma = \begin{pmatrix}
x & y & z & w \\
a & b & c & d \\
\end{pmatrix}$$

\noindent with $x \neq c$.  Take a sentential valuation $v$ such that $v(\epsilon aa) =  f$ and the rest is $t$. Then $v((\rm A_{S1}))$ = $f$ and $v(\sigma((\rm A_{S2})))$ = $t$.  

(Case 2) Let 
$$\sigma = \begin{pmatrix}
c & y & z & d \\
a & b & c & d \\
\end{pmatrix}.$$
Take a sentential valuation $v$ such that $v(\epsilon aa) =  v(\epsilon ad) = f$ and the rest is $t$.  Then $v((\rm A_{S1}))$ = $f$ and $v(\sigma((\rm A_{S2})))$ = $t$.  

(Case 3) Let 
$$\sigma = \begin{pmatrix}
c & y & z & b \\
a & b & c & d \\
\end{pmatrix}.$$
\noindent (Subcase 3.1) Let 
$$\sigma = \begin{pmatrix}
c & a & d & b \\
a & b & c & d \\
\end{pmatrix}.$$
In this case we see 

$\sigma((\rm A_{S2}))$ $= \sigma (\epsilon ab \wedge \epsilon cd . \supset . \epsilon cc \wedge (\epsilon bc \supset . \epsilon ad \wedge \epsilon ba))$

$\qquad \qquad = \enspace$ $\epsilon cd \wedge \epsilon ac . \supset . \epsilon aa \wedge (\epsilon da \supset . \epsilon bc \wedge \epsilon db)$.

Take a sentential valuation $v$ such that $v(\epsilon aa) =  v(\epsilon ac) = f$ and the rest is $t$. Then $v((\rm A_{S1}))$ = $f$ and $v(\sigma((\rm A_{S2})))$ = $t$. 
\medskip

\noindent (Subcase 3.2) Let 
$$\sigma = \begin{pmatrix}
c & d & a & b \\
a & b & c & d \\
\end{pmatrix}.$$
In this case we have 

$\sigma((\rm A_{S2}))$ $= \sigma (\epsilon ab \wedge \epsilon cd . \supset . \epsilon cc \wedge (\epsilon bc \supset . \epsilon ad \wedge \epsilon ba))$,

$\qquad \qquad = \enspace$ $\epsilon cd \wedge \epsilon ab . \supset . \epsilon aa \wedge (\epsilon da \supset . \epsilon cb \wedge \epsilon dc)$.
\medskip

In order to obtain $v((\rm A_{S1}))$ = $t$ and $v(\sigma((\rm A_{S2})))$ = $f$, we may have the following strategy:
\smallskip

$v(\epsilon ab \wedge \epsilon cd) = t$ and $v(\epsilon aa) = t$ for the common parts of $v((\rm A_{S1}))$ and $v(\sigma((\rm A_{S2})))$, 
\smallskip

$v(\epsilon bc \supset . \epsilon ad \wedge \epsilon ba) = t$ for $v((\rm A_{S1}))$, 
\smallskip

$v(\epsilon da \supset . \epsilon cb \wedge \epsilon dc) = f$ for $v(\sigma((\rm A_{S2}))$.

\noindent So, take a sentential valuation $v$ such that $v(\epsilon cb) = f$ (or $v(\epsilon dc) = f$) and the rest is $t$. Then we obtain $v((\rm A_{S1}))$ = $t$ and $v(\sigma((\rm A_{S2})))$ = $f$. 

(Case 4) Let 
$$\sigma = \begin{pmatrix}
c& x & y & a \\
a & b & c & d \\
\end{pmatrix}.$$
Take a sentential valuation $v$ such that $v(\epsilon aa) =  v(\epsilon da) =  v(\epsilon db) = v(\epsilon dc) = f$ and the rest is $t$. Then $v((\rm A_{S1}))$ = $f$ and $v(\sigma((\rm A_{S2})))$ = $t$

We shall prove (6). Recall 

\bigskip 

($\rm A_{S1}$) $\enspace$ $\epsilon ab \wedge \epsilon cd . \supset . \epsilon aa \wedge (\epsilon bc \supset . \epsilon ad \wedge \epsilon ba)$,

\bigskip   

($\rm A_{S3N}$) $\enspace$ $\epsilon ab \supset .\epsilon aa \wedge (\epsilon bc \supset . \epsilon bb \wedge (\epsilon cd \supset . \epsilon ad \wedge \epsilon ba))$. 

\bigskip 

We also have two cases for that.

(Case 1)  Let
$$\sigma = \begin{pmatrix}
x & y & z & w \\
a & b & c & d \\
\end{pmatrix}$$
with $y \neq a$.  Take a sentential valuation $v$ such that $v(\epsilon bb) = v(\epsilon cc) = v(\epsilon dd) = f$ and the rest is $t$. Then $v((\rm A_{S1}))$ = $t$ and $v(\sigma((\rm A_{S3N})))$ = $f$.  

(Case 2) Let 
$$\sigma = \begin{pmatrix}
x& a & z & w \\
a & b & c & d \\
\end{pmatrix}.$$
Take a sentential valuation $v$ such that $v(\epsilon bb) = f$ and the rest is $t$. Then $v((\rm A_{S1}))$ = $t$ and $v(\sigma((\rm A_{S3N})))$ = $f$.

We shall prove (9). Recall 

\bigskip 

($\rm A_{S2}$) $\enspace$ $\epsilon ab \wedge \epsilon cd . \supset . \epsilon cc \wedge (\epsilon bc \supset . \epsilon ad \wedge \epsilon ba)$,

\bigskip  

($\rm A_{S3N}$) $\enspace$ $\epsilon ab \supset .\epsilon aa \wedge (\epsilon bc \supset . \epsilon bb \wedge (\epsilon cd \supset . \epsilon ad \wedge \epsilon ba))$. 

\bigskip

We also have two cases for that.  

(Case 1)  Let
$$\sigma = \begin{pmatrix}
x & y & z & w \\
a & b & c & d \\
\end{pmatrix}$$
\noindent with $x \neq a$.  Take a sentential valuation $v$ such that $v(\epsilon aa) = v(\epsilon bb) = v(\epsilon dd) = f$ and the rest is $t$. Then $v((\rm A_{S2}))$ = $t$ and $v(\sigma((\rm A_{S3N})))$ = $f$.  

(Case 2) Let 
$$\sigma = \begin{pmatrix}
a & y & z & w \\
a & b & c & d \\
\end{pmatrix}.$$
Take a sentential valuation $v$ such that $v(\epsilon aa) = f$ and the rest is $t$. Then $v((\rm A_{S2}))$ = $f$ and $v(\sigma((\rm A_{S3N})))$ = $t$. 
\medskip


We shall prove (3d).  Recall

\bigskip 

($\mathrm{A_{M8}}$) $\enspace$ $\epsilon ab \wedge \epsilon cd . \supset . \epsilon aa \wedge \epsilon cc \wedge (\epsilon bc \supset . \epsilon ad \wedge \epsilon ba)$,

\bigskip

($\rm A_{S3Nd}$) $\enspace$ $\epsilon ab \supset .\epsilon aa \wedge (\epsilon bc \wedge \epsilon cd . \supset . \epsilon ad \wedge \epsilon ba)$. 

\bigskip

We have two cases for that.  

(Case 1)  Let
$$\sigma = \begin{pmatrix}
x & y & z & w \\
a & b & c & d \\
\end{pmatrix}$$
\noindent with $z \neq a$.  Take a sentential valuation $v$ such that $v(\epsilon cc) = f$ and the rest is $t$. Then $v((\rm A_{M8}))$ = $f$ and $v(\sigma((\rm A_{S3Nd})))$ = $t$.  

(Case 2) Let 
$$\sigma = \begin{pmatrix}
x & y & a & w \\
a & b & c & d \\
\end{pmatrix}.$$
Take a sentential valuation $v$ such that $v(\epsilon aa) = f$ and the rest is $t$. Then $v((\rm A_{M8}))$ = $f$ and $v(\sigma((\rm A_{S3Nd})))$ = $t$. 
\medskip

We shall prove (6d).  Recall

\bigskip 

($\rm A_{S1}$) $\enspace$ $\epsilon ab \wedge \epsilon cd . \supset . \epsilon aa \wedge (\epsilon bc \supset . \epsilon ad \wedge \epsilon ba)$,

\bigskip 

($\rm A_{S3Nd}$) $\enspace$ $\epsilon ab \supset .\epsilon aa \wedge (\epsilon bc \wedge \epsilon cd . \supset . \epsilon ad \wedge \epsilon ba)$. 

\bigskip

We have two cases for that.  

(Case 1)  Let
$$\sigma = \begin{pmatrix}
x & y & z & w \\
a & b & c & d \\
\end{pmatrix}$$
\noindent with $x \neq a$.  Take a sentential valuation $v$ such that $v(\epsilon aa) = f$ and the rest is $t$. Then $v((\rm A_{S1}))$ = $f$ and $v(\sigma((\rm A_{S3Nd})))$ = $t$.  

(Case 2) Let 
$$\sigma = \begin{pmatrix}
a & y & z & w \\
a & b & c & d \\
\end{pmatrix}.$$
Take a sentential valuation $v$ such that $v(\epsilon aa) = v(\epsilon cd) = f$ and the rest is $t$. Then $v((\rm A_{S1}))$ = $t$ and $v(\sigma((\rm A_{S3Nd})))$ = $f$. 
\medskip 


We shall prove (9d).  Recall 

\bigskip 

($\rm A_{S2}$) $\enspace$ $\epsilon ab \wedge \epsilon cd . \supset . \epsilon cc \wedge (\epsilon bc \supset . \epsilon ad \wedge \epsilon ba)$,

\bigskip

($\rm A_{S3Nd}$) $\enspace$ $\epsilon ab \supset .\epsilon aa \wedge (\epsilon bc \wedge \epsilon cd . \supset . \epsilon ad \wedge \epsilon ba)$. 

\bigskip

We have two cases for that.  

(Case 1)  Let
$$\sigma = \begin{pmatrix}
x & y & z & w \\
a & b & c & d \\
\end{pmatrix}$$
\noindent with $z \neq a$.  Take a sentential valuation $v$ such that $v(\epsilon cc) = f$ and the rest is $t$. Then $v((\rm A_{S2}))$ = $f$ and $v(\sigma((\rm A_{S3Nd})))$ = $t$.  

(Case 2) Let 
$$\sigma = \begin{pmatrix}
x & y & a & w \\
a & b & c & d \\
\end{pmatrix}.$$
Take a sentential valuation $v$ such that $v(\epsilon ab) = v(\epsilon cd) = f$ and the rest is $t$. Then $v((\rm A_{S2}))$ = $t$ and $v(\sigma((\rm A_{S3Nd})))$ = $f$. 
\medskip 

We shall prove (12d).  Recall 

\bigskip  

($\rm A_{S3N}$) $\enspace$ $\epsilon ab \supset .\epsilon aa \wedge (\epsilon bc \supset . \epsilon bb \wedge (\epsilon cd \supset . \epsilon ad \wedge \epsilon ba))$. 

\bigskip

($\rm A_{S3Nd}$) $\enspace$ $\epsilon ab \supset .\epsilon aa \wedge (\epsilon bc \wedge \epsilon cd . \supset . \epsilon ad \wedge \epsilon ba)$. 

\bigskip

We have two cases for that.  

(Case 1)  Let
$$\sigma = \begin{pmatrix}
x & y & z & w \\
a & b & c & d \\
\end{pmatrix}$$
\noindent with $y \neq a$.  Take a sentential valuation $v$ such that $v(\epsilon bb) = f$ and the rest is $t$. Then $v((\rm A_{S3N}))$ = $f$ and $v(\sigma((\rm A_{S3Nd})))$ = $t$.  

(Case 2) Let 
$$\sigma = \begin{pmatrix}
x & a & z & w \\
a & b & c & d \\
\end{pmatrix}.$$
Take a sentential valuation $v$ such that $v(\epsilon ab) = v(\epsilon aa) = f$ and the rest is $t$. Then $v((\rm A_{S3N}))$ = $f$ and $v(\sigma((\rm A_{S3Nd})))$ = $t$. $\Box$ 
\medskip


The idea for ($\rm A_{M8}$), ($\rm A_{S1}$), ($\rm A_{S2}$), ($\rm A_{S3N}$) and ($\rm A_{S3Nd}$) is that $\epsilon bc$ is taken as an intermediary to connect $\epsilon ab$ and $\epsilon cd$ for the transitivity (Ax2).

\section{Conjecture 1. $A_k$-type single axiom schemata}

As candidates of nontrivial axiom schemata (on 2025-2-2), we shall present certain axiom schemata to be considered as follows.

\bigskip 

($\rm A_{k1}$) $\enspace$ $\epsilon ab \supset . \epsilon aa \wedge (\epsilon bb \wedge \epsilon bc . \supset . \epsilon ac \wedge \epsilon ba)$. 

\bigskip 

($\rm A_{k2}$) $\enspace$ $\epsilon ab \supset . \epsilon aa \wedge (\epsilon cc \wedge \epsilon bc . \supset . \epsilon ac \wedge \epsilon cb)$.

\bigskip 

($\rm A_{k3}$) $\enspace$ $\epsilon ab \supset . \epsilon aa \wedge (\epsilon cd \wedge \epsilon bc . \supset . \epsilon ac \wedge \epsilon cb)$.

\section{Conjecture 2. Additional single axiom schemata, Part I}

As candidates of nontrivial axiom schemata (on 2025-2-2), we shall give certain axiom schemata to be considered as follows.
\bigskip 

($\rm A_{ad1}$) $\enspace$ $\epsilon ab \wedge \epsilon bb . \supset . \epsilon aa \wedge \epsilon ba \wedge (\epsilon bc \supset \epsilon ac)$.

\bigskip 

($\rm A_{ad2}$) $\enspace$ $\epsilon ab \supset . \epsilon aa \wedge (\epsilon bc \supset \epsilon ac) \wedge (\epsilon bb \supset \epsilon ba)$. 

\bigskip

($\rm A_{ad6}$) $\enspace$ $\epsilon ab \wedge \epsilon bc . \supset . \epsilon aa \wedge \epsilon ba \wedge (\epsilon cd \supset \epsilon bd)$.

\bigskip

($\rm A_{ad6-2}$) $\enspace$ $\epsilon ab \wedge \epsilon bc . \supset . \epsilon bb \wedge \epsilon ba \wedge (\epsilon bd \supset \epsilon ad)$.

\bigskip

($\rm A_{ad7}$) $\enspace$ $\epsilon ab \wedge \epsilon bc . \supset . \epsilon aa \wedge \epsilon ba \wedge (\epsilon cd \supset \epsilon ad)$.

\bigskip

($\rm A_{ad7-2}$) $\enspace$ $\epsilon ab \wedge \epsilon bc . \supset . \epsilon bb \wedge \epsilon ba \wedge (\epsilon cd \supset \epsilon ad)$.

\bigskip 

($\rm A_{ad8}$) $\enspace$ $\epsilon ab \wedge \epsilon bc . \supset . \epsilon aa \wedge \epsilon bb \wedge \epsilon ac \wedge \epsilon ba$.

\section{Conjecture 3. Additional single axiom schemata, Part II}

As candidates of nontrivial axiom schemata (on 2025-2-2), we shall give six axiom schemata to be considered as follows.
\bigskip 

($\rm A_{S1ex1}$) $\enspace$ $\epsilon ab \wedge \epsilon cd . \supset . \epsilon aa \wedge (\epsilon bc \supset . \epsilon bd \wedge \epsilon ba)$.
\bigskip

($\rm A_{S1ex2}$) $\enspace$ $\epsilon ab \wedge \epsilon cd . \supset . \epsilon aa \wedge (\epsilon bc \supset . \epsilon bd \wedge \epsilon cb)$.
\bigskip

($\rm A_{S1ex3}$) $\enspace$ $\epsilon ab \wedge \epsilon cd . \supset . \epsilon aa \wedge (\epsilon bc \supset . \epsilon ac \wedge \epsilon cb)$.
\bigskip

($\rm A_{S2ex1}$) $\enspace$ $\epsilon ab \wedge \epsilon cd . \supset . \epsilon cc \wedge (\epsilon bc \supset . \epsilon bd \wedge \epsilon ba)$.
\bigskip

($\rm A_{S1ex2}$) $\enspace$ $\epsilon ab \wedge \epsilon cd . \supset . \epsilon cc \wedge (\epsilon bc \supset . \epsilon bd \wedge \epsilon cb)$.
\bigskip

($\rm A_{S1ex3}$) $\enspace$ $\epsilon ab \wedge \epsilon cd . \supset . \epsilon cc \wedge (\epsilon bc \supset . \epsilon ac \wedge \epsilon cb)$.
\bigskip

\section{Summary and Remarks}
On March 8, 1995, was found the following nontrivial single axiom schema characteristic of Le\'{s}niewski-Ishimoto's propositional ontology $\mathbf{L_1}$  (Inou\'{e} \cite{ inoue16}).
$$(\mathrm{A_{M8})} \enspace \epsilon ab \wedge \epsilon cd . \supset . \epsilon aa \wedge \epsilon cc \wedge (\epsilon bc \supset . \epsilon ad \wedge \epsilon ba).$$
 (For recent work on Le\'{s}niewski's system, see e.g., Indrzejczak \cite{Indrzejczak2022}, Urbaniak \cite{urbaniak-book} and Inou\'{e} \cite{inoue-Blass, inoue2021}.) ) The original paper did not provide a definition and proof of the nontriviality of $(\mathrm{A_{M8}})$. We will present a definitive definition and proof of it here, along with an update on the progress concerning this axiom schema since 1995. For this purpose, we introduce two novel criteria, nontriviality and quasi-nontriviality, to distinguish between two axiom schemata.

In proving quasi-nontriviality, the subformulas of the form $\epsilon aa$ play essential roles in demonstrating quasi-nontriviality, in principle.

As main results, we will present simplified axiom schemata $(\mathrm{A_{S1}})$, $(\mathrm{A_{S2}})$, $(\mathrm{A_{S3N}})$, and $(\mathrm{A_{S3Nd}})$ based on $(\mathrm{A_{M8}})$, and discuss their nontriviality and quasi-nontriviality.

This study opens the way to consider the use of computers to obtain further advances in this field, as is the trend in algebra (refer to Kunen \cite{Kunen1996a, Kunen1996b}, Britten et al. \cite{Britten}, Phillips and Vojt\v{e}chovsk\'{y} \cite{PhillipsVojt}, etc.).

We believe that Le\'{s}niewski's system will increasingly become important in studying the foundations of mathematics, particularly concerning the existence of mathematical objects and the future extension of mathematics itself.

This research also contributes to the humanization of mathematics, a theme explored further in Inou\'{e} \cite{inoue-human}.

The first author of this paper thinks that there is still a possibility to add some single axiom schemata for $\bf L_1$.

\bigskip 

\noindent Takao Inou\'{e}

\noindent Faculty of Informatics

\noindent Yamato University

\noindent Katayama-cho 2-5-1, Suita, Osaka, 564-0082, Japan

\noindent inoue.takao@yamato-u.ac.jp
 
\noindent (Personal) takaoapple@gmail.com (I prefer my personal mail)

\bigskip

\noindent Tadayoshi Miwa

\noindent Library, The University of Tokyo Library

\noindent Hongo 7-3-1, Bunkyo-ku, Tokyo 113-0033, Japan

\noindent miwa.tadayoshi@mail.u-tokyo.ac.jp

\end{document}